\documentclass[12pt]{article}
\usepackage{amssymb}
\usepackage[all]{xy}
\font\tenfrak=eufm10
\font\sevenfrak=eufm7
\font\fivefrak=eufm5
        \newfam\frakfam  
                \textfont\frakfam=\tenfrak
\font\tenDDl=msbm10  
\font\sevenDDl=msbm7 
\font\fiveDDl=msbm5 
        \newfam\DDlfam \def\DDl{\fam\DDlfam\tenDDl} 
                \textfont\DDlfam=\tenDDl 
\scriptfont\DDlfam=\sevenDDl
        \scriptscriptfont\DDlfam=\fiveDDl
\scriptfont\frakfam=\sevenfrak  
        \scriptscriptfont\frakfam=\fivefrak

 \def\Fhd#1#2{\smash{\mathop{\hbox to 14mm{\rightarrowfill}}
\limits^{\scriptstyle#1}_{\scriptstyle#2}}}

\def\Fhg#1#2{\smash{\mathop{\hbox to 14mm{\leftarrowfill}}
\limits^{\scriptstyle#1}_{\scriptstyle#2}}}

 \def\fhd#1#2{\smash{\mathop{\hbox to 8mm{\rightarrowfill}}
\limits^{\scriptstyle#1}_{\scriptstyle#2}}}

\def\fhg#1#2{\smash{\mathop{\hbox to 8mm{\leftarrowfill}}
\limits^{\scriptstyle#1}_{\scriptstyle#2}}}

\def\fhnull#1#2{\smash{\mathop{\hbox to 0mm{
}}
\limits^{\scriptstyle#1}_{\scriptstyle#2}}}

\def\fhUp#1#2#3{\smash{\mathop{\hbox to 8mm{
$#2$}}
\limits^{\scriptstyle#1}_{\scriptstyle#3}}}

\def\diagram#1{\def\normalbaselines{\baselineskip=0pt
\lineskip=0pt\lineskiplimit=0pt}   \matrix{#1}}

\long\def\nodo#1{}
\def\genfd{{\bf k}}
\def\mat22#1#2#3#4{ \left(\begin{array}{cc} #1 & #2\\ #3 & #4
\end{array}\right) } 
\newcommand{\id}{{\rm id}}
\def\dj{d\kern-.30em\raise1.25ex\vbox{\hrule width .3em height .03em}}
\def\Dj{D\rlap{\kern-.70em\raise0.75ex
    \vbox{\hrule width .3em height .03em}}}
\def\comdG{{\bf G}}
\def\comdg{G}

\long\def\nodo#1{}
\long\def\detail#1{} 
\long\def\notforme#1{} 
\newif\ifzsnoco
\newcounter{point}

\def\lnamedef#1{\expandafter\edef\csname lpt#1 \endcsname}
\def\lnameuse#1{\expandafter\csname lpt#1 \endcsname}

\def\End{{\rm End}\,}

\begin{document}
\begin{center}
{\Large Cyclic structures for simplicial objects from comonads} \\
{\it preliminary notes}\\
\footnotesize{
{\sc Zoran \v{S}koda}, email: {\tt zskoda@irb.hr}  \vskip .03in
{Theoretical Physics Division}\\
{Institute Rudjer Bo\v{s}kovi\'{c}, P.O.Box 180}
} 
\end{center}

\vskip .2in
{\footnotesize 
The simplicial endofunctor induced
by a comonad in some category may underly
a cyclic object in its category of endofunctors.
The cyclic symmetry is then given by a sequence of
natural transformations. We write down 
the commutation relations the first cyclic operator
has to satisfy with the data of the comonad.
If we add a version of quantum Yang Baxter relation
and another relation we actually get a sufficient
condition for constructing a sequence of 
higher cyclic operators in a canonical fashion.
A degenerate case of this construction comes
from so-called trivial symmetry of an additive comonad.

We also consider weaker versions for paracyclic objects
as well as some connections to the subject of distributive laws.
}

\vskip .02in

{\footnotesize {\sc AMS classification:} {\bf 19D55}}

\vskip .2in

Comonadic homology and monadic ('triple') cohomology are among 
the standard unifying frameworks in homological algebra. 
Despite great effort, I could not find 
{\it any} references giving instances of cyclic (co)homology 
derived completely in that same framework 
(that is, including the cyclic operator as well).  
There is maybe one exception: operads are
often treated in monadic language, and 
one enriched version, called {\it cyclic operads}, is designed 
in part as a gadget 
to do the cyclic homology~(\cite{cyclicoper, Menichi:BV}). 

We perform here a pretty straightforward general nonsense exercise
to provide a (nontautological) additional structure $t$
on an arbitrary comonad $\comdG$ in ${\cal A}$ generating the 
cyclic operators $t_n$ in the sense of
{\sc Connes}~(\cite{Connes:book,Loday:cyclicbk, Weibel})
on the simplicial objects associated to $\comdG$. 
We also relate our data to some other functorial 
data, e.g. a class of distributive laws.
Our procedure is just one of at least several conceivable ways 
one could utilize (co)monads to obtain cyclic-type constructions. 
This kind of production of cyclic objects is functorial. 
Although this was intended (I won't dwell
on motivating picture as it is largely still conjectural), 
at the present moment this is a drawback: 
most known kinds of cyclic cohomology do
not have coefficients, hence should not fit into our framework.
I hope that the proposed framework is transparent enough to enable us
to spot new classes of examples, remedying the problem.  

{\bf Prerequisites.} Given categories ${\cal A,B,C}$,
functors $f_1,f_2,g_1,g_2$ and natural transformations $F,G$ as in the diagram
$$\diagram{
      {\cal A}\,\,\, \fhUp{\fhd{f_1}{}}{\Uparrow \!\!\!F}{\fhd{}{f_2}}
\,\,\,{\cal B}\,\,\, \fhUp{\fhd{g_1}{}}{\Uparrow \!\!\!G}{\fhd{}{g_2}}
\,\,\,{\cal C},
}$$
define the natural transformation $G\star F : g_2 \circ f_2
\Rightarrow g_1 \circ f_1$ by
\[
 ( G\star F )_A := G_{f_1(A)} \circ g_2(F_A) = g_1(F_A) \circ G_{f_2(A)}
:  g_2(f_2(A)) \rightarrow g_1 (f_1(A)).\]
$(F,G) \mapsto F\star G$ is called the {\bf Godement product}
('horizontal composition').
It is associative for triples for which $F\star(G\star H)$ is defined.

A {\bf (co)monad} in category ${\cal A}$ 
is a (co)monoid in the monoidal category ${\rm End}\,{\cal A}$
of endofunctors in ${\cal A}$~(\cite{Semtriples,MacLane,Weibel}).
The monoidal product of endofunctors is their composition and
of natural transformations the Godement product.  
Equivalently, a comonad is a triple 
$\comdG = (\comdg,\delta,\epsilon)$ where $\comdg : {\cal A}\to {\cal A}$
is an endofunctor in  ${\cal A}$ ('underlying endofunctor of $\comdG$'), 
and the 'comultiplication' $\delta : \comdG \Rightarrow \comdG^2$ 
and the 'counit' $\epsilon : \comdG \Rightarrow 1_{\cal A}$ 
are natural transformations of functors, 
such that for every object $M$ in ${\cal A}$
the coassociativity axiom
$\delta_{GM}\circ \delta_M = G(\delta_M) \circ \delta_M$ 
and the counit axiom 
$\epsilon_{GM} \circ \delta_M = G(\epsilon_M) \circ \delta_M = 1_{GM}$  
hold.

Let ${\bf \Delta}$ be the 'cosimplicial' category:
its objects are nonnegative integers viewed as
finite ordered sets ${\bf n} := \{ 0 < 1 <  \ldots < n\}$ and 
its morphisms are nondecreasing monotone functions.
Given a category ${\cal A}$, denote by ${\rm Sim}{\cal A}$
the category of {\bf simplicial objects} in ${\cal A}$, i.e. functors
$F : {\bf \Delta}^{\rm op}\rightarrow {\cal A}$.
Represent $F$ in ${\rm Sim}{\cal A}$ as a sequence 
$F_n := F({\bf n})$ of objects, 
together with the face maps $\partial_i^n : A_n \rightarrow A_{n-1}$
and degeneracy maps $\sigma_i^n : F_n \rightarrow F_{n+1}$
for $i \in {\bf n}$ satisfying
the familiar simplicial identities~(\cite{MacLane,Weibel}). 
Notation $F_\bullet$ for these data is standard.

To any comonad $\comdG$ in ${\cal A}$ one associates
the sequence $\comdG_\bullet$ of endofunctors
${\DDl Z}_{\geq 0} \ni n \mapsto \comdG_n := \comdg^{n+1} := 
\comdg \circ \comdg \circ \ldots \circ \comdg$,
together with natural transformations 
$\partial_i^n : \comdg^i \epsilon \comdg^{n-i} : 
\comdg^{n+1} \rightarrow \comdg^n$
and $\sigma_i^n : \comdg^i \delta \comdg^{n-i} : 
\comdg^{n+1} \rightarrow \comdg^{n+2}$,
satisfying the simplicial identities. 
Hence any comonad $\comdG$ canonically induces a 
simplicial endofunctor\index{simplicial endofunctor},
i.e. a functor $\comdG_\bullet : \Delta^{\rm op} 
\rightarrow {\rm End}{\cal A}$, or equivalently, a functor 
$\comdG_\bullet : {\cal A} \rightarrow {\rm Sim}{\cal A}$.
The counit $\epsilon$ of the comonad $\comdG$ satisfies
$\epsilon \circ \partial_0^1 = \epsilon \circ \partial_1^1$, 
what means that 
$\epsilon : \comdG_\bullet \rightarrow {\rm Id}_{\cal A}$ 
is in fact an augmented simplicial endofunctor.

A ${\DDl Z}$-cyclic (synonym: paracyclic) 
object in ${\cal A}$ is a simplicial object $F_\bullet$ together with
a sequence of isomorphisms $t_n : F_n \rightarrow F_n$, $n\geq 1$,
such that
\begin{equation}\label{eq:cycrel}\begin{array}{cc}
\partial_i t_n = t_{n-1} \partial_{i-1},\,\, i > 0, &
\sigma_i t_n = t_{n+1} \sigma_{i-1},\,\, i > 0, \\
\partial_0 t_n = \partial_n, & \sigma_0 t_n = t_{n+1}^2 \sigma_n.
\end{array}\end{equation}
A ${\DDl Z}$-cocyclic (paracocyclic) object in
${\cal A}$ is a ${\DDl Z}$-cyclic object in ${\cal A}^{\rm op}$.
${\DDl Z}$-(co)cyclic object is (co)cyclic if, 
in addition, $t_n^{n+1} = 1$ 
(\cite{Connes:book, Loday:cyclicbk, Weibel}).
Equivalently, the category $\Delta^{\rm op}$ may be upgraded to
the {\bf cyclic category} 
${\cal C}$ of Connes~(\cite{Connes:book,Loday:cyclicbk}). 
It is the universal category containing $\Delta^{\rm op}$
as a nonfull subcategory, identical on objects, and having 
minimal set of additional morphisms containing 
a sequence of ``cyclic'' morphisms $\tau_n : {\bf n}\to {\bf n}$ such that
any simplicial object $F_\bullet$ in any category $\cal A$ is a cyclic
if the operators $t_n := F(\tau_n)$ are declared cyclic. 
\vskip .02in

{\bf Bottom relations.} 
Let now $\comdG = (\comdg,\delta,\epsilon)$ 
be a comonad on a category ${\cal A}$ and 
$t : \comdg \comdg \Rightarrow \comdg \comdg$ a natural transformation. 

For every object $M$ in ${\cal A}$ we require
\begin{equation}\label{eq:tcomonadcomrelM}\begin{array}{rclcl}
 \comdg(\epsilon_M) t_M  & = & \epsilon_{\comdg M} 
& : & \comdg^2 M \rightarrow \comdg M, \\
 \epsilon_{\comdg M} t_M  & = & \comdg(\epsilon_M) 
& : & \comdg^2 M \rightarrow \comdg M,\\
 \comdg(\delta_M) t_M & = & t_{\comdg M} \comdg(t_M) \delta_{\comdg M} 
& : & \comdg^2 M \rightarrow \comdg^3 M, \\
 \delta_{\comdg M} t_M & = &  
t_{\comdg M} \comdg(t_M) t_{\comdg M} \comdg(t_M) \comdg (\delta_M)
& : & \comdg^2 M \rightarrow \comdg^3 M,
\end{array}\end{equation}
or, in a more schematic form of natural transformations,
\begin{equation}\label{eq:tcomonadcomrel}\begin{array}{rclcl}
 \comdg(\epsilon) t & = & \epsilon_{\comdg} 
\\
 \epsilon_{\comdg } t  & = & \comdg(\epsilon) 
\\
 \comdg(\delta) t & = & t_{\comdg } \comdg(t) \delta_{\comdg} \\
 \delta_{\comdg } t & = &  
[t_{\comdg} \comdg(t)]^2 \comdg (\delta)
\end{array}\end{equation}
Notice that if $t_M \circ t_M = {\rm Id}_{\comdg\comdg M}$ then
the first two identities in~(\ref{eq:tcomonadcomrelM}) 
are equivalent (composing one by $t_M$ from the right). These relations
are the bottom-part of the relations required for the cyclic symmetry.
Now we show that with few additional properties they are sufficient, 
as well.

{\bf Theorem 1.} {\it 
Let  $\comdG = (\comdg,\delta,\epsilon)$ 
be a comonad on a category ${\cal A}$,
where $\epsilon$ is the counit and $\delta$ the coproduct. Let
$\comdG_\bullet$ be the associated simplicial endofunctor.
Suppose an invertible natural transformation
 $t : \comdg \comdg \Rightarrow \comdg \comdg$ satisfies 
~(\ref{eq:tcomonadcomrel}), and the quantum Yang Baxter equation (QYBE)
\begin{equation}\label{eq:QYBEendof}
 G(t_{M})\circ t_{GM} \circ G(t_{M})
 = t_{GM}\circ G(t_{M}) \circ t_{GM}, 
\,\,\,\,\forall M\in {\rm Ob}\,{\cal A}.
\end{equation}
If we also assume the relation
\begin{equation}\label{eq:tcubedeqtcubed}
t_M \circ t_M\circ t_M\circ \delta_M = t_M \circ \delta_M,
 \,\,\, \,\,\,\forall M\in {\rm Ob}\,{\cal A},
\end{equation}
then setting
\begin{equation}\label{eq:constrtn}
 t_{nM} : = 
 t_{G^{n-1}M} \circ G(t_{G^{n-2}M})\circ \ldots  \circ G^{n-1}(t_M), 
\end{equation}
defines paracyclic operators $t_n$ 
on the augmented simplicial endofunctor
$\comdG_\bullet \stackrel{\epsilon}\Rightarrow {\rm Id}_{\cal A}$
making it into an augmented paracyclic object in $\End {\cal A}$.
}

{\bf Remark.} 1. We do not claim that every paracyclic 
or even every cyclic operator on $\comdG_\bullet$ is of that kind.

2. Practical (weaker) form of the conclusion: 
for any fixed object $M$ in ${\cal A}$, the pair
$({\bf G}_\bullet M \stackrel{\epsilon_M}\rightarrow M, t_{*M})$
is an augmented paracyclic object in ${\cal A}$.

{\bf Observation.} 
\begin{equation}\label{eq:tn}
t_{n+1\,M} = t_{n\,GM} \circ G^n(t_M).
\end{equation}

{\it Proof of the theorem.} Substituting 
$\partial_i = G^i(\epsilon_{G^{n-i}})$,  $\sigma_i = G^i(\delta_{G^{n-i}})$
and~(\ref{eq:constrtn}) in~(\ref{eq:cycrel}) we obtain
\begin{itemize}
\item $(A_{n,i})$  for  $i = 1,\ldots, n$,
\[ G^i (\epsilon_{G^{n-i}M}) t_{G^{n-1}M}\ldots G^{n-1}(t_M)    
=  t_{G^{n-2}M}\ldots G^{n-2}(t_M) 
G^{i-1}(\epsilon_{G^{n-i+1}M}),\]
\item $(B_{n,i})$ for  $i = 1,\ldots, n$,
\[ G^i (\delta_{G^{n-i}M}) t_{G^{n-1}M}\ldots G^{n-1}(t_M)
=  t_{G^{n}M}\ldots G^{n}(t_M) G^{i-1}(\delta_{G^{n-i+1}M}),\]
\item $(C_n)$ $\epsilon_{G^n M} t_{G^{n}M}\ldots G^{n-1}(t_M) 
=  G^n(\epsilon_M)$,
\item $(D_n)$ $\delta_{G^n M} t_{G^{n}M}\ldots G^{n-1}(t_M) 
=  [t_{G^{n}M}\ldots G^{n}(t_M)]^2 G^n(\delta_M)$.
\end{itemize}
Of course, the factors involving $G^{n-2}$ in our notation 
appear only if $n>1$.

Basis of induction: 
$(A_{1,1})$ is the first and $(C_1)$ the second formula in
~(\ref{eq:tcomonadcomrel}) while $(B_1)$ is the third and 
$(D_1)$ the bottom formula there.

The rest will follow by inductive calculations 
(as usual RHS= right-hand side etc.).
Cases A,B,C are very simple:

a) $(A_{n,i}) \Rightarrow (A_{n+1,i+1})$ Act by $G$ on both sides of
equation  $(A_{n-i})$, and compose both sides 
from the left by $t_{G^{n-1}M}$.
Then use the naturality formula 
$t_{G^{n-1}M} G^{i+1} (\epsilon_{G^{n-i}M}) 
= G^{i+1} (\epsilon_{G^{n-i}M}) t_{G^n M}$ on the left-hand side.

b) $(A_{n,i}) \Rightarrow (A_{n+1,i})$ 
Write down $(A_{n,i})$ on $GM$ instead of $M$:
\[
G^i (\epsilon_{G^{n+1-i}M}) t_{G^{n}M}\ldots G^{n-1}(t_{GM})    = 
 t_{G^{n-1}M}\ldots G^{n-2}(t_{GM}) 
G^{i-1}(\epsilon_{G^{n-i+2}M}),
\]
and compose from right both sides by $G^{n}(t_M)$. At RHS use
the naturality formula
$G^{i-1}(\epsilon_{G^{n-i+2}M}) G^{n}(t_M) = 
G^{n-1}(t_M)G^{i-1}(\epsilon_{G^{n-i+2}M})$, which holds for $n-i>0$.

c) $(B_{n,i}) \Rightarrow (B_{n+1,i+1})$ -- analogously to a).

d)  $(B_{n,i}) \Rightarrow (B_{n+1,i})$ -- Write down $(B_{n,i})$ on $GM$
instead of $M$ and compose from right both sides by $G^{n}(t_M)$.
At RHS use the naturality formula
$G^{i-1}(\delta_{G^{n-i+2}M}) G^{n}(t_M) = 
G^{n+1}(t_M)G^{i-1}(\delta_{G^{n-i+2}M})$ for $n-i > 0$.
 
e) Assume~$(C_n)$. Then
$$\xymatrix{ G^{n} G G \ar[r]_{G^{n}(t)}  \ar@/^/[rr]^{t_{n+1}}
\,\,\,\,\ar[rd]_{G^{n}G(\epsilon)} 
& \,G^{n} G G \ar[r]_{t_{n,G}}\, \ar[d]^{\! G^{n}(\epsilon_G)} 
& \,\,\,\,G G^{n} G \ar[dl]^{\,\,\,\epsilon_{G^{n}G}} \\ 
& G^{n+1} &
}$$
is commutative: the left triangle 
by $G(\epsilon) = \epsilon_G(t)$ and the functoriality of $G^n$,
and the right triangle is commutative by the inductive hypothesis 
(composed by $G$). The external triangle is then~$(C_{n+1})$.

f) Finally, the case $(D_n) \Rightarrow (D_{n+1})$ 
is somewhat more elaborate, and this is the point 
where we also need QYBE. For this, we first notice that $(D_1)$, together
with QYBE, implies the following calculation:
\[ \begin{array}{lcl}
G^n(\delta_{GM}) G^n(t_M) &=& G^n(\delta_{GM} t_M) \\
&\stackrel{(D_1)}{=} & G^n\left[t_{GM}G(t_M) 
t_{GM}G(t_M) G(\delta_M)\right]\\ 
&=&  G^n(t_{GM}) G^{n+1}(t_M)G^n(t_{GM}) G^{n+1}(t_M)G^{n+1}(\delta_M)\\
&\stackrel{{\rm QYBE}}{=}& G^{n+1}(t_M)G^n(t_{GM}) 
G^{n+1}(t_M) G^{n+1}(t_M) G^{n+1}(\delta_M) \\
 &=& G^{n+1}(t_M) G^n(t_{GM})G^{n+1}(t_M t_M \delta_M).
\end{array}\]
Thus, 
\begin{equation}\label{eq:Gnpqybe}\begin{array}{rcl}
G^n(t_{GM}) G^n(\delta_{GM}) G^n(t_M) &=&
G^n(t_{GM})G^{n+1}(t_M) G^n(t_{GM}) G^{n+1}(t_M^2\delta_M)\\
\mbox{= (use QYBE)} &=& G^{n+1}(t_M)G^n(t_{GM}) G^{n+1}(t_M^3 \delta_M)\\
\mbox{= (use~(\ref{eq:tcubedeqtcubed}))}
&=& G^{n+1}(t_M)G^n(t_{GM}) G^{n+1}(t_M) G^{n+1}(\delta_M)
\end{array}\end{equation}
We write down the $(D_n)$ for $GM$ instead of $M$, and then 
compose both sides by $G^{n}(t_M)$ from the right:
$$\begin{array}{l}
\delta_{G^{n+1} M} t_{{n}\,GM} G^{n} (t_M) 
= [t_{n+1\,GM}]^2 G^n(\delta_{GM}) G^{n}(t_M),
\\
\delta_{G^{n+1} M} t_{{n+1}\,M} = t_{n+1\,GM} t_{n\,GM} 
G^n(t_{GM}) G^n(\delta_{GM}) G^{n}(t_M).
\end{array}$$
Now we substitute the identity from~(\ref{eq:Gnpqybe}) to get
\[
\delta_{G^{n+1} M} t_{{n+1}\,M} = t_{{n+1}\,GM} t_{{n}\,G^2 M}
G^{n+1}(t_M)G^n(t_{GM})G^{n+1}(t_M)G^{n+1}(\delta_M).
\]
and then we notice that naturality and the definition of
$t_{n}$ imply that $t_{n\,GM}$ commutes with
$G^{n+1}(t_M)$. Hence
\[\begin{array}{lcl}
\delta_{G^{n+1} M} t_{{n}\,M} &=& 
t_{n+1\,GM} G^{n+1}(t_{M})t_{n\,GM} 
G^n(t_{GM}) G^{n+1}(t_M)G^{n+1}(\delta_M)\\
&=&  [t_{n+2\,M}]^2 G^{n+1}(\delta_M).
\end{array}\]
{\bf Theorem 2.} {\it Assume in addition that $t^2 = 1$.
Then $t_n^{n+1} = 1$ as well, hence
$\comdG_\bullet \stackrel{\epsilon}\Rightarrow {\rm Id}_{\cal A}$
is an augmented {\bf cyclic} object in $\End {\cal A}$.
}

{\bf Remark.} More generally, under the conditions of Theorem 1,
given a fixed object $M$ in ${\cal A}$, 
for the augmented paracyclic object 
$({\bf G}_\bullet M \stackrel{\epsilon_M}\rightarrow M, t_{*M})$
to be a cyclic object in ${\cal A}$ it is sufficient that 
$t^2_{G^k M} = {\rm Id}_M$ for all $k\geq 0$. 

{\it Proof of Theorem 2}. 
Theorem 1 being proved it remains to verify $t_n^{n+1} = 1$. 
However this is standard. 
Namely, QYBE and the naturality of $t$ imply that transformations 
$\alpha^i_n = G^i(t_{G^{n-i-1}}) : G^n  \to G^n$ for $0\leq i \leq n$
are the standard generators (braids) of 
the braid group $B_{n+1}$ on $n+1$ letters in a representation into
by natural autoequivalences of functor $G^n$
\[
   B_n \to {\rm Aut}(G^n).
\]
If $t^2 = 1$ then $(\alpha^i_n)^2 = [G^i(t_{G^{n-i-1}})]^2 = 1$ as well. 
It is then a standard result,
that this representation factors through
the symmetric group on $n+1$ letters.  
In this representation, $t_n$ by its definition~(\ref{eq:constrtn})
equals to $\alpha^0_n \alpha^1_n  \ldots \alpha^n_n$,
what is easily recognized as 
the image of a standard cycle in the symmetric group of order $n+1$.
Alternatively, one could use inductive calculations 
to show  $(\alpha^0_k \alpha^1_k  \ldots \alpha^n_k)^{n+1}=1$
whenever $k\geq n$. A convenient intermediate step is to show 
$(\alpha^0_k  \alpha^1_k  \ldots \alpha^{n-1}_k\alpha^n_k)^{n-1} =
(\alpha^0_k  \alpha^1_k  \ldots \alpha^{n-1}_k)^{n-1} 
\alpha^{n}_k\ldots \alpha^1_k$.

{\bf Definition.} 
{\it Given a comonad $\comdG$ in ${\cal A}$, a natural transformation
$t : \comdg\comdg\rightarrow \comdg\comdg$ is }
\begin{itemize}
\item (\cite{grandis:symmsimp}) {\bf a symmetry} {\it of $\comdG$ 
if it sastisfies the QYBE~(\ref{eq:QYBEendof}) and also}
\begin{equation}\label{eq:symmG}
t^2 = 1_{GG},\,\,\,\,t\delta = \delta,\,\,
\,\,\,\epsilon_G t = G(\epsilon),\,\,\,
\delta_G t = G(t) t_G G(\delta)\end{equation}
{\it In that case, the pair $(\comdG,t)$ is called a {\bf symmetric comonad.}}

\item {\bf a strong braiding} {\it of $\comdG$ 
if it sastisfies the QYBE~(\ref{eq:QYBEendof}) and also}
\begin{equation}\label{eq:ZsymmG}
G(\epsilon) t = \epsilon_G,\,\,\,\epsilon_G t = G(\epsilon),
\,\,\,t_G G(t) \delta_G = G(\delta) t,\,\,\,
\delta_G t = G(t) t_G G(\delta)\end{equation}
\end{itemize}

{\bf Lemma.} {\it Every symmetry of a comonad is a strong braiding.} 

{\it Proof.} It is immediate that~(\ref{eq:symmG}) imply
$t_G G(t) \delta_G = G(\delta) t$ and $G(\epsilon) t = \epsilon_G$.
Relation $t^2 = 1_{GG}$ and QYBE imply together 
$G(t) t_G  =  [t_G G(t)]^2$. 
\vskip .02in

{\bf Observation.} Every symmetric comonad $(\comdG,t)$ satisfies
the conditions of Theorem 2, and 
hence it gives a rule for producing certain cyclic objects in ${\cal A}$.
Conversely, every $t$ satisfying the conditions of the Theorem 2,
and satisfying $t^2 = 1$ is a symmetry. 

{\bf Remark.} A strong braiding on a comonad does not 
imply the conditions of Theorem 1. Namely,
$[t_{\comdg M} \comdg(t_M)]^2 \comdg (\delta_M)\neq
\comdg(t_M)t_{\comdg M}\comdg (\delta_M)$ in general.

{\bf Remark.} Symmetric comonads are related to
{\it symmetric simplicial sets}~(\cite{grandis:symmsimp}). 

{\bf Proposition 1.} (cf.~\cite{Menini:MSRI})
{\it If ${\cal A}$ is an {\em additive category}, 
and comonad $\comdG$ respects
its additive structure, then 
$\tau = \tau^G : GG \to GG$
defined by
$\tau_M := \delta_M \epsilon_{GM} + \delta_M G(\epsilon_M) - 1_{GGM}$
is a symmetry of $\comdG$.

$\tau^G$ is called the {\bf trivial symmetry} or simply 
{\bf the symmetry} of $\comdG$.
}

{\it Proof.} $t\delta = \delta$ and $\epsilon_G t = G(\epsilon)$
are immediate by the counit axioms.
Other relations are left to the reader -- 
they follow by calculations involving many summands 
(particularly the QYBE). Use the naturality of $\delta$
and $\epsilon$ and axioms for comonad when collecting and comparing
the summands. Q.E.D.
\vskip .17in

Now we compute $t_n = \tau_{G^{n-1}}\circ G(t_{n-1})$
where $\tau$ is the trivial symmetry. 
As this is a combinatorial problem,
we will introduce some helpful notation. 

Define a small strict monoidal category ${\cal P}$ as follows. 
Objects of ${\cal P}$ are nonnegative integers 
(written with square brackets $[n]$). 
Morphisms ${\cal P}([m],[n])$ are $m$-tuples $(k_1,\ldots, k_m)$
such that $\sum_{i=1}^m k_i = n$. The tensor product on objects is the 
addition $([m],[n]) \mapsto [m+n]$. The tensor product of morphisms
is the concatenation. We denote it either by $\star$ or simply concatenate. 
If ${\bf r} = (r_1,\ldots, r_m), {\bf k} = (k_1,\ldots, k_n)$ 
are tuples which are composable (i.e. $\sum_{i=1}^n k_i = m$)
then we define ${\bf r}\circ {\bf k}$ as follows. 
Represent all components $r_i$ ($k_j$) by trees of height one 
with $r_i$ ($k_j$ respectively) leaves. 
Attach $r_1$ to the left most leaf in ${\bf k}$ 
(of course we jump over trees corresponding to zeros as they are leafless).
Then attach $r_2$ to the next leaf to the right from that leaf (on the same
or on the next tree) and so on. In other words, attach $r_i$ to the
$i$-th leaf from the left in ${\bf k}$
(this leaf is at the $(k_1 + \ldots + k_{i-1} +1)$-th tree 
from the left in the graphical presentation of ${\bf k}$).
Then in the resulting depth $\leq 2$ trees remove the internal nodes
(in particular at the places where the zero is attached the whole leg
disappears after that proces). The pictures with variants of such
short trees are useful for extension
of similar calculations to more complicated setups.

Equivalently, in nongraphical presentation, the composition can be
described as follows. ${\bf r}\circ {\bf k}$ has as many nodes as ${\bf k}$. 
The number at $i$-th place is the sum of $k_i$ components
taken successively from ${\bf r}$ starting at position 
$1 + \sum_{p = 0}^{i-1} k_p$ and proceeding 
taking components to the right from there.

{\bf Example.} $(2031041)\circ (021202)=(0,2+0,3,1+0,0,4+1)=(023105)$.

{\bf Proposition 2.} {\it ${\cal P}$ is a small strict monoidal category. }

In particular, the composition and 
the tensor product are strictly associative, 
and we have the interchange law
$(r\circ r')\star(s\circ s') = (r\star s)\circ (r'\star s')$ 
holds for morphisms. The proof is an exercise.

For any $i\geq 0$ introduce the natural transformation
$\Delta_i : G\to G^i$ as follows:
$\Delta_0 := \epsilon$, $\Delta_1 := \id_G$, $\Delta_2 := \delta$,
and inductively $\Delta_{k+1} := \delta_{G^{k-1}}\circ\Delta_k$ for $k> 2$.
We define a strict monoidal functor $\Delta : {\cal P}\to {\rm End}\,{\cal A}$
as follows. On objects $\Delta([n]):= G^i$. On morphisms,
$\Delta(k) = \Delta_k$ for $k \in {\cal P}([1],[k])$ (in fact $(k)$).
As any morphism $s$ in ${\cal P}([m],[n])$ is an $m$-tuple, hence a
tensor product (concatenation) of $m$ $1$-tuples in ${\cal P}([1],[s_i])$,
there is at most one extension of these formulas 
respecting {\it monoidal} structures, namely  
$\Delta({\bf k})=\Delta_{\bf k}:= \Delta_{k_1}\star\ldots\star\Delta_{k_n}$ 
where ${\bf k} = (k_1,\ldots,k_n)$. 
The reader may verify that $\Delta$ is in fact a functor i.e. 
$\Delta(s\circ s') = \Delta(s)\circ\Delta(s')$ where $\circ$ 
on the left-hand side is the composition of endofunctors. 
A crucial part of that statement is to use the coassociativity 
of $\delta$ and that the property that $\epsilon$. 
is a counit translates to the fact $01\circ 2 = 10\circ 2 = 1$.
One also has $k'\circ k = k+k'-1$ for $k>0$, and in particular
$2\circ k = k+1$.
We conclude:

{\bf Proposition 3.} {\it $\Delta : {\cal P}\to {\rm End}\,{\cal A}$ 
is a monoidal functor.}

\vskip .1in
Notice that set 
$\coprod {\cal P}([m],[n]) =\coprod_{n\geq 0} {\DDl Z}^n_{\geq 0}$
is an ordered monoid with respect to the concatenation
and lexicographic ordering. $\Delta$ also encodes
the map, also denoted by $\Delta : \coprod_{n\geq 0} {\DDl Z}^n_{\geq 0}\to
{\rm End}\,{\cal A}$ which evaluate on $n$-tuples 
the same as the functor $\Delta$ evaluates on
the corresponding morphisms in ${\cal P}$.

Most important part in the following are the endomorphism
sets $S_n := {\cal P}([n],[n])$ which are monoids with respect to
the composition. This are sets of $n$-tuples
of nonnegative integers which moreover add up to $n$.
The restriction 
$\Delta = \Delta| : \coprod_n S_n \to \cup_n {\rm Nat}(G^n,G^n)$ 
is a map of monoids.

The composition in ${\cal P}$ 
restricted to $S_n\times S_n$ takes values in $S_n$. 
Hence it is an associative binary operation
and it extends to a bilinear operation on ${\DDl Z}S_n$.

Since ${\cal A}$ is an abelian, hence additive category, 
${\rm Nat}(G^n,G^n)$ is an abelian group in particular.
Therefore $\Delta|_{S_n}$ extends to a unique homomorphism
from the free abelian group ${\DDl Z}S_n$ with basis $S_n$ to
${\rm Nat}(G^n,G^n)$, also denoted by $\Delta$.

Let ${\rm NAlt}_n$ be the subset of ${\DDl Z}S_n$ consisting of
all $x$ of the alternating sign form $x = s_1 - s_2 + s_3 - \ldots \pm s_k$
where $k\geq 1$, all $s_i \in S_n$ and $s_1 < s_2 < \ldots < s_k$
with respect to the lexicographic ordering. For example,
writing $0004$ for $(0,0,0,4)$ etc. the element
$0004 - 0040 + 1210 - 1300$ is in ${\rm NAlt_n}$.

We also form $S^\tau_n \subset S_n$ inductively, along with
a linear map $\alpha : {\DDl Z}S^\tau_n\to {\DDl Z}\tilde{S}_{n}$ 
where $\tilde{S}$ is the
set of all sequences of $n-1$ characters from the alphabet $\{a,b,c\}$.
We set $S^\tau_2 := \{ 02,11,20 \}$, $\alpha(02):=a$, $\alpha(11):=b$,
$\alpha(20):=c$. If $S^\tau_{n-1}$ is defined, then

1. $(0,s_2,\ldots,s_n)\in S^\tau_n$ iff 
$\left(s_2 > 1 \mbox{ and } (s_2-1, \ldots,s_n)\in S^\tau_{n-1}\right)$.
\newline In that case, $\alpha(0,s_2,\ldots,s_n) = a\alpha(s_2,\ldots,s_n)$.

2. $(1,s_2,\ldots,s_n)\in S^\tau_n$ iff
$(s_2,\ldots,s_n)\in S^\tau_{n-1}$.
\newline In that case, $\alpha(0,s_2,\ldots,s_n) = b\alpha(s_2,\ldots,s_n)$.

3. $(2,s_2,\ldots,s_n)\in S^\tau_n$ iff $\left(s_2=s_3 =\ldots= s_{n-1} = 1
\mbox{ and }s_n = 0\right)$. 
\newline In that case, $\alpha(2,s_2,\ldots,s_n) = cc\ldots c$ ($n-1$ times).
\vskip .1in

Finally, ${\rm NAlt}^\tau_n := {\rm NAlt}_n \cap {\DDl Z}S^\tau_n$. 

A $k$-tuple of sequences $s_1,\ldots,s_k \in S_n$ form a {\bf twin $k$-tuple} 
if they become identical after subtracting 1 from the first nonzero member
of each sequence. For example, $0030$ and $1020$ form a twin pair, because 
they both become $0020$ after performing this operation.

{\bf Definition.} $x\in {\DDl Z}S^\tau_n$ 
{\bf represents} $u\in {\rm Nat}(G^n,G^n)$ if $\Delta(x) = u$.
Then $x$ is called below the representation and
$\alpha(x)$ the  $abc$-representation of $u$.

{\bf Theorem 3.} {\it 1. $S^\tau_n$ has exactly $2^{n-1}+1$ elements.

2. a sequence in $\tilde{S}$ of $n-1$ characters belongs to
the image of $\alpha$ iff characters $a$ and $b$ never appear after $c$.

3. Let ${\bf m} = {\bf m}_1 - {\bf m}_2 +\ldots + {\bf m}_{2^{n-1}+1}
\in{\rm NAlt}^\tau_n$, ${\bf m}_1<\ldots<{\bf m}_{2^{n-1}+1}$
be the unique element of ${\rm NAlt}^\tau_n$ of maximal length $2^{n-1}+1$. 
Then $\Delta({\bf m})=t_{n-1}$.

4. If ${\bf m}_k$ is the $k$-th summand in ${\bf m}$, beginning in $0$
then ${\bf m}_{k+2^{n-2}}$ is its twin 
which comes with an opposite sign in ${\bf m}$.

5. The sign of ${\bf m}_k$ is positive iff character $b$
appears even number of times in $\alpha({\bf m}_k)$.
}

{\bf Example.} $t_4 = \Delta( 00005-00041+00050-00302+00311-00320
+00410-02003+02021-02030+02102-02111+02120-02210+03110
-10004+10031-10040+10202-10310+11003-11021+11030-11102+11111
-11120+11210-12110+21110).$

{\it Proof of the theorem.} (sketch) We prove all assertions simultaneously
by induction. The case $n=2$: $t_1 = \Delta(02-11+20)
= \epsilon\star\delta - \id + \delta\star\epsilon$.
Together with the definition of $\Delta$ this also implies that
$t_{G^{n-1}} = \Delta(02(1^{n-1})-11(1^{n-1})+20(1^{n-1}))$,
where $1^{n-1}=11\ldots 1$ ($n-1$ times).
 
Suppose all 5 assertions hold for $n-1$. Then 
$G(t_{n-1}) = \sum_{i = 1}^{2^{n-2}+1} \Delta(1\star{\bf m}_i)$, i.e.
at the level of ${\DDl Z}S^\tau_n$ we attach $1$ from the left to
each summand in ${\bf m}$ (the element representing $t_{n-1}$). 
For example, $G(t_1) = \Delta(102-111+120)$. 
The reader should play enough to get comfortable calculating with
the composition in ${\cal P}$.
In any case, composing $02(1^{n-1})$ from the left, 
decrements the leftmost nonzero number in
a sequence and increments the first next nonzero number to the right in it;
and $20(1^{n-1})$ does other way around. As the leftmost number in the
representation of $G(t_{n-1})$ is $1$, it becomes $0$ after that
composition. Hence composing by
$02(1^{n-1})$ results in twins of the original summands representing 
$G(t_{n-1})$, with opposite sign. For $20(1^{n-1})$, 
the increments and decrements are switched and we end with a sequence
starting with $2$ and the rest is truncated in a way to erase
the difference between the twins; so all such terms, pairwise cancel
except for the twinless $211\ldots10$. As a result we end without
any other sequences headed by $2$; and all new sequences 
not headed by $2$ are paired in twins.
Verifying the rules for signs and remaining requirements is straightforward.
\vskip .1in

Let $\genfd$ be a commutative ring, $p,r\in \genfd$,
${\cal A}$ a $\genfd$-linear category, 
and ${\bf G}$ a $\genfd$-linear comonad in ${\cal A}$.
Following~\cite{BrzNich} we define
$\theta := p\epsilon\star\delta - p\id + r\delta\star\epsilon$.
Then $\theta$ generalizes the trivial symmetry (just set $p=r=1$), 
it still satisfies the QYBE, but it is not a symmetry of ${\bf G}$,
and it is not even a strong braiding (as the reader
should check all 4 equations in~(\ref{eq:ZsymmG}) fail in general). 
As it satisfies the QYBE, the braid cycle 
$t_n$ may still be of interest to compute. 
The basic combinatorics from the $\tau$-case actually passes through! 

 If $t_1 = \theta$ then $t^\theta_n = t_n$ is
described as follows. Suppose $s_i \in S^\tau_n$ corresponds to
the sequence $a_1 a_2 \ldots a_k c^{n-k}$ where $a_i \in \{a,b\}$
and $k = k(i)$ are determined by $s_i$ as explained above.
Define a map ${\rm norm}_\theta : {\DDl Z}S^\tau_n\to {\genfd}S^\tau_n$ by
setting $s_i \mapsto p^{k(i)} r^{n-k(i)}$ and extending additively.
Let $S^\theta_n := {\rm norm}_\theta(S^\tau_n)$,
${\rm NAlt}_n^\theta :={\rm norm}_\theta({\rm NAlt}_n^\tau)$
and ${\bf m}^\theta = \theta({\bf m})$ where 
${\bf m}\in {\rm NAlt}_n^\tau$ is the longest element as above. 
$\Delta$ extends $\genfd$-linearly to a map 
$\Delta :\genfd S_n\to {\rm Nat}(G^n,G^n)$ and 
this map restricts down to $\genfd S^\theta_n$ and ${\rm NAlt}^\theta_n$.
\vskip .03in

{\bf Proposition 4.} {\it $t^\theta_n = \Delta({\bf m}^\theta)$.
}
\vskip .01in

{\it Proof.} The proof of the Theorem 1, passes through, but one
has to keep track of coefficients. It is easy to see that the
coefficients of each summand are easily
tracked in terms of the $abc$-representation. 
The key observation is that those twin pairs 
which cancel in the inductive step of Theorem 1
actually do come with exactly the same coefficients (but different
signs, as before) so they again cancel after aplying $20(1)^{n-1}$
in the next step. 
\vskip .02in

{\bf Definition.}~(\cite{Barr:composite, Beck:distr})
A {\bf distributive law} from a comonad
${\bf G} = (G,\delta^G, \epsilon^G)$ to a comonad
${\bf F} = (F,\delta^F, \epsilon^F)$ is a natural
transformation $l : F\circ G\to G\circ F$ such that 
\[
\begin{array}{cc}G(\epsilon^F)\circ l = (\epsilon^F)_{G}, &
G(\delta^F) \circ l = l_F \circ F(l) \circ (\delta^F)_G,\\
(\epsilon^G)_{F} \circ l = F(\epsilon^G), &
(\delta^G)_F \circ l = G(l)\circ l_G \circ F(\delta^G).\end{array} 
\]
Any distributive law $l$ induces a {\bf composite comonad}
\[ 
{\bf F}\circ_l {\bf G} = 
(F\circ G,\delta^{l,F\circ G},\epsilon^{l,F\circ G}),
\]
where  the coproduct $\delta^{l,F\circ G}$ equals the composition 
\[ FG \stackrel{F(\delta^G)}\longrightarrow
 FGG \stackrel{(\delta^F)_{GG}}\longrightarrow FFGG
\stackrel{F(l_{G})}\longrightarrow FGFG,\]
and the counit
$\epsilon^{l,F\circ G}$ equals the composition 
$FG \stackrel{F(\epsilon^G)}\longrightarrow 
F\stackrel{\epsilon^F}\longrightarrow 1$.

{\bf Theorem 4.} {\it 
Suppose $t_1 = t$ is a strong braiding on ${\bf G}$.
Define inductively $t_{n+1 M} := t_{n GM} \circ G^n(t_M)$. 
Then $t_n : G^n \circ G \Rightarrow G \circ G^n$ 
form a sequence of distributive laws from ${\bf G}$
to ${\bf G}^n = (G^n, \delta^{(n-1)}, \epsilon^{(n-1)})$, 
where  the comultiplication $\delta^{(n)}$ is inductively defined by
\begin{equation}\label{eq:coprn}
\delta^{(n)} = G^{n-1}(t_{n-1, G})\delta^{(n-1)}_{G^2} G^{n-1}(\delta)
= G^{n-1}(t_{n-1, G}) G^{2n-2}(\delta)\delta^{(n-1)}_G,
\end{equation}
and $\delta^{(1)} = \delta$. The counit 
is $\epsilon^{(n)} = \epsilon^{(n-1)}_G \circ G^{n-1}(\epsilon) =
\epsilon \circ G(\epsilon)\circ\cdots\circ G^{n-1}(\epsilon)$.

In formulas, the following holds for $n\geq 1$:
\begin{eqnarray}\label{eq:hidist1}
G(\delta^{(n)}) t_n = t_{n, G^n} G^n(t_n) \delta^{(n)}_G\\
\label{eq:hidist2}
\delta_{G^n} t_n = G(t_n) t_{n,G} G^n(\delta)\\
\label{eq:hidist3}
(\epsilon^{(n)})_G = G(\epsilon^{(n)}) t_n\\
\label{eq:hidist4}
G^n(\epsilon) = \epsilon_{G^n} t_n
\end{eqnarray}
Actually, (\ref{eq:hidist1}) is a case $l = n$ 
of a more general identity:
\begin{equation}\label{eq:hidist1gen}
G^{n-l+1}(\delta^{(l)}) t_n = t_{n+l} G^{n-l}(\delta^{(l)}_G),
\,\,\,\,1\leq l \leq n
\end{equation}
}
The converse does not hold: a family of distributive laws as above
does not need to be coming from the first plus QYBE.

{\bf Lemma.} (a) {\it Let $t^{(1)} = t : GG \Rightarrow GG$ 
be any natural transformation and let $t^{(n)}$ be inductively defined by
$t^{(n+1)} := (t^{(n)})_G G^n(t)$, $n \geq 1$, then 
\begin{equation}\label{eq:lema}
t^{(n+i)} = t^{(n)}_{G^i} G^n (t^{(i)})
\,\,\mbox{ for all }\,\,1\leq i, 1\leq n.
\end{equation}
}
(b) {\it Under the assumptions from (a) and the QYBE~(\ref{eq:QYBEendof}),
\begin{equation}\label{eq:lemb}
G(t^{(n-1)}_{G}) t^{(n+1)}= t^{(n+1)} t^{(n-1)}_{G^2}
\,\,\mbox{ holds for }\,\,n\geq 2.
\end{equation}
More generally,
\begin{equation}\label{eq:lembp}
G^p(t^{(l-1)}_{G}) t^{(p+l)}= t^{(p+l)} G^{p-1}(t^{(l-1)}_{G^2})
\,\,\mbox{ holds for }\,\,l\geq 2, p\geq 1.
\end{equation}
}
(c) {\it If in addition $[t^{(1)}]^2 = \id$, 
then $[t^{(n+1)}]^2 = G(t^{(n)}) t^{(n)}_G$ holds for $n\geq 1$.
}
\vskip .03in

{\it Proof.} (a) For $i=1$ this is the definition of $t_{n+1}$.
Suppose the lemma holds for some $i$. Then
$t_{n+(i+1)} = t_{n+i,G} G^{n+i}(t) = t_{n, G^{i+1}} G^n(t_{i,G}) G^{n+i}(t)
= t_{n, G^{i+1}} G^n(t_{i,G}G^i(t)) = t_{n,G^{i+1}} G^n(t_{i+1})$.

(b) These are (equivalent to) simple identities in the braid group.
For completeness, we give a direct proof. 

Base of induction: for $n = 2$ we have the identity
\begin{equation}\label{eq:lemb2}
G(t_G) t^{(3)} = t^{(3)} t_{G^2},
\end{equation}
which follows by applying the QYBE:
$G(t_G) t^{(3)} = G(t_G)t_{G^2}G(t_G)G^2(t) = 
t_{G^2}G(t_G)t_{G^2}G^2(t)=t_{G^2}G(t_G)G^2(t)t_{G^2}=t^{(3)}t_{G^2}$.

$$\begin{array}{lcl}
G(t^{(n-1)}_G) t^{(n+1)} &=& G(t^{(n-2)}_{G^2})G^{n-1}(t_{G^2})
t^{(n-2)}_{G^3}G^{n-2}(t_{G^2})G^{n-1}(t_G)G^n(t)
\\&=&G(t^{(n-2)}_{G^2})t^{(n-2)}_{G^3}G^{n-1}(t_{G^2})
G^{n-2}(t_{G^2})G^{n-1}(t_G)G^n(t)
\\&\stackrel{{\rm QYBE}}=&G(t^{(n-2)}_{G^2})t^{(n-2)}_{G^3}
G^{n-2}(t_{G^2})G^{n-1}(t_G)G^{n-2}(t_{G^2})G^n(t)
\\&=&G(t^{(n-2)}_{G^2})t^{(n-2)}_{G^3}
G^{n-2}(t_{G^2})G^{n-1}(t_G)G^n(t)G^{n-2}(t_{G^2})
\\&=&G(t^{(n-2)}_{G^2})t^{(n)}_G G^n(t)G^{n-2}(t_{G^2})
\\&=&\left[G(t^{(n-2)}_G)t^{(n)}\right]_G G^n(t)G^{n-2}(t_{G^2})
\\\mbox{(induction on $l$)}&=&\left[G(t^{(n-l-1)}_G)t^{(n-l+1)}\right]_{G^l} 
G^n(t)G^{n-2}(t_{G^2}), \,\,\,1\leq l \leq n-2,
\\\left(l = n-2\right)&=&\left[G(t_G)t^{(3)}\right]_{G^{n-2}}
G^3(t_{G^{n-3}})\cdots G^n(t) G(t_{G^{n-1}})\cdots G^{n-2}(t_{G^2})
\\&\stackrel{(\ref{eq:lemb2})}=& t^{(3)}_{G^{n-2}} t_{G^n} G^3(t_{G^{n-3}})
\cdots G^n(t) \cdot G(t_{G^{n-1}})\cdots G^{n-2}(t_{G^2})
\\&=& t^{(3)}_{G^{n-2}} G^3(t_{G^{n-3}})
\cdots G^n(t) \cdot t_{G^n} G(t_{G^{n-1}})\cdots G^{n-2}(t_{G^2})
\\&=& t^{(n+1)} t^{(n-1)}_{G^2}.
\end{array}$$

$$\begin{array}{lcl}
G^p(t^{(l-1)}_G) t^{(p+l)} 
&=&G^p(t^{(l-1)})t^{(p-1)}_{G^{l+1}} G^{p-1}(t^{(l+1)})
\\&=&t^{(p-1)}_{G^{l+1}}G^p(t^{(l-1)}_G)G^{p-1}(t^{(l+1)})
\\&=&t^{(p-1)}_{G^{l+1}}G^{p-1}(G(t^{(l-1)}_G)t^{(l+1)})
\\&\stackrel{(\ref{eq:lemb})}=&
t^{(p-1)}_{G^{l+1}}G^{p-1}(t^{(l+1)}t^{(p-1)}_{G^2}))
\\&=&[t^{(p-1)}_{G^{l+1}}G^p(t^{(l+1)})]G^{p-1}(t^{(p-1)}_{G^2}))
\\&\stackrel{(\ref{eq:lema})}=&t^{(p+l)}G^{p-1}(t^{(l-1)}_{G^2}).
\end{array}$$

(c) is a simple identity in the symmetric group $\Sigma(n+2)$.

{\it Proof of Theorem 4}. To warm up, we start with $n = 2$ case of
eq.~(\ref{eq:hidist1}). 
{\it In calculations, we will write $t^{(n)}$ for $t_n$.}

\begin{equation}\label{eq:delta2}
\delta^{(2)} = G(t_G)\delta_{G^2} G(\delta) = G(t_G) G^2(\delta)\delta_G
\end{equation}

$$\begin{array}{lcl}
G(\delta^{(2)}) t^{(2)} &\stackrel{(\ref{eq:delta2})}=& 
G^2 (t_G) G^3(\delta) G(\delta_G) t_G G(t)\\
&=& G^2(t_G) G^3(\delta) (G(\delta) t)_G G(t) \\
&=& G^2(t_G) G^3(\delta) t_{G^2} G(t_G) \delta_{G^2} G(t) \\
&=& G^2(t_G) t_{G^3} G^3(\delta)G(t_G) G^2(t) \delta_{G^2}\\
&=& t_{G^3} G^2(t_G) G(t_{G^2}) G^3(\delta) G^2(t) \delta_{G^2}\\
&=& t_{G^3} G^2(t_G) G(t_{G^2}) G^2(G(\delta)t) \delta_{G^2}\\
&=& t_{G^3} G^2(t_G) G(t_{G^2}) G^2(t_G G(t) \delta_G) \delta_{G^2}\\
&=& t_{G^3} G^2(t_G) G(t_{G^2}) G^2(t_G) G^3(t) G^2(\delta_G) \delta_{G^2}\\
&\stackrel{{\rm QYBE}}=&
t_{G^3} G(t_{G^2}) G^2(t_{G}) G(t_{G^2}) G^3(t) G^2(\delta_G) \delta_{G^2}\\
&=&
t_{G^3} G(t_{G^2}) G^2(t_{G}) G^3(t) G(t_{G^2}) G^2(\delta_G) \delta_{G^2}\\
&=&
t^{(4)} (G(t_{G}) G(\delta_G) \delta_{G})_G\\
&\stackrel{(\ref{eq:delta2})}=& t^{(4)} \delta^{(2)}_G.
\end{array}$$

The equation~(\ref{eq:hidist1}) is just the
special case of~(\ref{eq:hidist1gen}) when $l=n$. 
Now we prove~(\ref{eq:hidist1gen}). 

$$\begin{array}{lcl}
G^{n-l+1}(\delta^{(l)}) t^{(n)} &=& G^{n-l-1}
(G^{l-1}(t^{(l-1)}_G)G^{2l-2}(\delta)\delta^{(l-1)}_G)t^{(n-1)}_G G^{n-1}(t)
\\&=&G^n(t^{(l-1)}_G)G^{n-l+1}(\delta)
G^{n-l+1}(\delta^{(l-1)}_G)t^{(n-1)}_G G^{n-1}(t) 
\\&=&G^n(t^{(l-1)}_G)G^{n-l+1}(\delta)
\left[G^{n-1-(l-1)+1}(\delta^{(l-1)})t^{(n-1)}\right]_G G^{n-1}(t)
\\&=&G^n(t^{(l-1)}_G)G^{n-l+1}(\delta)
t^{(n+l-2)}_G G^{n-l}(\delta^{(l-1)}_{G^2})G^{n-1}(t)
\\&=&G^n(t^{(l-1)}_G)t^{(n+l-2)}_{G^2} G^{n+l-1}(\delta)
G^{n+l-2}(t)G^{n-l}(\delta^{(l-1)}_{G^2})
\\&=&G^n(t^{(l-1)}_G)t^{(n+l-2)}_{G^2} G^{n+l-2}(G(\delta)t)
G^{n-l}(\delta^{(l-1)}_{G^2})
\\&=&G^n(t^{(l-1)}_G)t^{(n+l-2)}_{G^2} 
G^{n+l-2}(t^{(2)})G^{n+l-2}(\delta_G)G^{n-l}(\delta^{(l-1)}_{G^2})
\\&=&G^n(t^{(l-1)}_G)t^{(n+l)}G^{n+l-2}(\delta_G)G^{n-l}(\delta^{(l-1)}_{G^2})
\\&\stackrel{(\ref{eq:lembp})}=&
t^{(n+l)}G^{n-1}(t^{(l-1)}_{G^2})
G^{n+l-2}(\delta_G)G^{n-l}(\delta^{(l-1)}_{G^2})
\\&=&t^{(n+l)}G^{n-l}(G^{l-1}(t^{(l-1)}_G)G^{2l-2}(\delta)\delta_G^{(l-1)})_G
\\&=&t^{(n+l)}G^{n-l}(\delta^{(l)}_G).
\end{array}$$

If $t^2 = \id$ then the lemma, part (c) holds and
~(\ref{eq:hidist2}) is identical to $(D_n)$ 
in the proof of Theorem 1, as well as the assumptions from Theorem 1
are also fullfilled. 
In general, eq.~(\ref{eq:hidist2}) is different from $(D_n)$, 
but we will prove it by induction, following very similar path as for $(D_n)$.
\begin{equation}\label{eq:Gnpqybe2}\begin{array}{rcl}
G^n(t_G) G^n(\delta_G) G^n(t) &=&
G^n(t_G)G^n(\delta_G t)\\&=&
G^n(t_{G})G^n\left[G(t)t_G G(\delta)\right]\\&=&
G^n(t_{G})G^{n+1}(t) G^n(t_G) G^{n+1}(\delta)\\
\mbox{= (use QYBE)} &=& G^{n+1}(t)G^n(t_G) G^{n+1}(t) G^{n+1}(\delta)\\
\mbox{= (use~(\ref{eq:tcubedeqtcubed}))}
&=& G^{n+1}(t)G^n(t_G) G^{n+1}(t) G^{n+1}(\delta)
\end{array}\end{equation}
We write down the~(\ref{eq:hidist2}) pushed by $G$ from the right, 
and then compose both sides by $G^{n}(t)$ also from the right:
$$\begin{array}{l}
\delta_{G^{n+1}} t^{(n)}_G G^{n} (t) 
= [t^{(n+1)}_G]^2 G^n(\delta_G) G^{n}(t),
\\
\delta_{G^{n+1}} t^{(n+1)} = t^{(n+1)}_G t^{(n)}_G  
G^n(t_G) G^n(\delta_G) G^{n}(t).
\end{array}$$
We substitute the identity proved in~(\ref{eq:Gnpqybe2}) to get
$$
\delta_{G^{n+1}} t^{(n+1)} = t^{(n+1)}_G t^{(n)}_{G^2}
G^{n+1}(t)G^n(t_G)G^{n+1}(t)G^{n+1}(\delta).
$$
and then we notice that naturality and the definition of
$t^{(n)}$ imply that $t^{(n)}_G$ commutes with $G^{n+1}(t)$. Hence
\[\begin{array}{lcl}
\delta_{G^{n+1}} t^{(n)} &=& 
t^{(n+1)}_G G^{n+1}(t)t^{(n)}_G 
G^n(t_G) G^{n+1}(t)G^{n+1}(\delta)\\
&=&  [t^{(n+2)}]^2 G^{n+1}(\delta).
\end{array}\]

Eq.~(\ref{eq:hidist3})  follows by induction, and naturality. 
Assume~(\ref{eq:hidist3}) holds for $n$. Then
$$\xymatrix{
G^n G G \ar[rr]^{G^n(t)} \ar[dr]_{G^n(\epsilon_G)}
&& G^n G G \ar[rr]^{t^G_n}\ar[dl]^{G^{n+1}(\epsilon)} 
&& G G^n G\ar[dl]^{G^{n+1}(\epsilon)}\\
& G^n G \ar[rr]_{t_n}\ar[dr]_{\epsilon^{(n)}_G} 
&& GG^n\ar[dl]^{G(\epsilon^{(n)})} & \\
&& G &&
}$$
is commutative (upper right triangle by $G(\epsilon)t=\epsilon_G$
and functoriality of $G^n$; the parallelogram by naturality of $t_n$,
and the lower triangle by the induction hypothesis). 
The external triangle may easily be recognized as the identity 
$(\epsilon^{(n+1)})_G = G(\epsilon^{(n+1)}) t_{n+1}$,
once $\epsilon^{(n+1)}$ is expanded as $\epsilon^{(n)}\circ G^n(\epsilon)$
and $t_{n+1}$ as $t_{n,G} \circ G^n(t)$.

Eq.~(\ref{eq:hidist4}) is identical to~$(C_{n+1})$ from the proof of
Theorem 1. Although the assumptions on $t$ differ slightly, the
proof stays unchanged.  Q.E.D.
\vskip .1in

Notice that in the special case $t^2 = \id$, 2 out of 4 
requirements for the distributive laws are the identities 
for the cyclic operator on the nose, but the other two differ.
\vskip .1in

Let $f : R\to S$ be an extension of noncommutative
unital rings. Consider the $S$-bimodule 
$S^{\otimes_R n} = S\otimes_R S\otimes_R \cdots \otimes_R S$
and let $\tau$ be the trivial symmetry on $S$. Given an $S$-bimodule map
$r : S\otimes_R S \to {}_S M_{S}$ let
$r_{ii+1} : S^{\otimes_R k}\to S^{\otimes_R i-1}\otimes_R 
M\otimes S^{\otimes_R k-i-1}$, $1\leq i \leq k$ 
be the $S$-bimodule map obtained by evaluating $r$ 
on $i$-th and $i+1$-st places and the identity 
on the rest of tensor factors.

{\bf Theorem 5.} {\it
The formula (written in {\sc Nuss}~\cite{Nuss}) 
\begin{equation}\label{eq:Nuss}\begin{array}{lcl}
\mu^{(n)} &=& (\mu_{12}\circ \mu_{34}\circ \ldots \circ \mu_{2n-1,2n})
\circ \tau_{2n-2,2n-1} \circ (\tau_{2n-4,2n-3}\circ \tau_{2n-3,2n-2})
\\&&\circ\, (\tau_{2n-4,2n-3}\circ \tau_{2n-3,2n-2}\circ \tau_{2n-2,2n-1})
\circ \ldots \\&&\circ\, (\tau_{45}\circ\ldots \circ \tau_{n+1,n+2})
\circ (\tau_{23}\circ\ldots\circ\tau_{nn+1}),
\end{array}\end{equation}
for the induced associative multiplication on $S^{\otimes_R n}$, 
is (a special case of) the dualization of our formula for
coproduct~(\ref{eq:coprn}) expanded by induction.
}

{\it Proof.} We iterate formula~(\ref{eq:coprn}) descending $n$:
$$\begin{array}{lcl}
\delta^{(n)} &=& G^{n-1}(t^{(n-1)}_G) G^{n-2}(t^{(n-2)}_{G^3})
\delta^{(n-2)}_{G^4}G^{n-2}(\delta_{G^2})G^{n-1}(\delta)\\
&=& G^{n-1} (t^{(n-1)}_G) G^{n-2}(t^{(n-2)}_{G^3})
\circ\ldots\\&&\,\circ\, G^{n-p}(t^{(n-p)}_{G^{2p-1}})\delta^{(n-p)}_{2p}
G^{n-p}(\delta_{G^{2p}})G^{n-p+1}(\delta_{G^{2p-2}})\circ\ldots\circ 
G^{n-1}(\delta)\\&=&G^{n-1} (t^{(n-1)}_G) G^{n-2}(t^{(n-2)}_{G^3})
\circ\ldots\\&&\,\circ\, G(t_{2n-3})\delta_{2n-2}
G(\delta_{G^{2n-4}})G^2(\delta_{G^{2n-6}})\circ\ldots\circ G^{n-1}(\delta)
\\&=&\prod_{j=1}^{n-1} G^{n-j}(t^{(n-j)}_{G^{2j-1}})
\prod_{p=0}^{n-1} G^p(\delta_{G^{2n-2p-2}})\end{array}$$
Now we need to dualize; the arrows and the 
composition will be hence backwards. Dualizing
$G^{n-j}(t^{(n-j)}_{G^{2j-1}}) = 
G^{n-j}(t_{G^{n+j-2}}) G(t_{G^{n+j-1}})\ldots G^{2n-2j-1}(t_{G^{2j-1}})$
will hence be
$(\tau_{2j,2j+1}\circ\tau_{2j+1,2j+2}\circ\ldots\circ\tau_{n+j-1,n+j})$.
Similarly, the dualization of $G^{p}(\delta_{G^{2n-2p-2}})$ is
$\mu_{2n-2p-1,2n-2p}$. 
Putting these together, in proper order, we get~(\ref{eq:Nuss}).
As distributive laws were designed in 1960s for exactly that kind of 
reason, we believe that the formula and our explanation for it,
must have been known to experts before.

{\bf Lemma.} {\it Let $\cal C, \cal A$ be any two categories,
and ${\cal G}:{\cal C}\to{\rm End}\,{\cal A}$, 
${\cal X}:{\cal C}\to{\cal A}$ functors, either both covariant
or both contravariant. Then the rule
$$\begin{array}{l}
{\cal Y}_C := {\cal Y}(C) 
:= {\cal G}(C)({\cal X}(C)), \,\,\,\,\,C \in {\rm Ob}\,{\cal C},\\
{\cal Y}_f :={\cal Y}(f) 
:= {\cal G}(C')({\cal X}(f))\circ({\cal G}(f))_{{\cal X}(C)}
= ({\cal G}(f))_{{\cal X}(C')}\circ{\cal G}(C)({\cal X}(f)),
\end{array}$$
for $f \in {\cal C}(C,C')$, 
defines a functor ${\cal Y} : {\cal C}\to{\cal A}$ of the same covariance.
}

{\it Proof.} This is basically the same trick 
which is involved in the definition of Godement's product. 
Given a chain $C\stackrel{f}\to C'\stackrel{f}\to C''$ in ${\cal C}$, 
the diagram
$$\xymatrix{
{\cal G}(C)({\cal X}(C))
\ar[d]^{{\cal G}(C)({\cal X}(f))}
\ar[rrd]^{{\cal Y}_f}
&&&&
\\{\cal G}(C)({\cal X}(C'))\ar[rr]^{{\cal G}(f)_{{\cal X}(C')}}
\ar[d]^{{\cal G}(C)({\cal X}(g))}
&&{\cal G}(C')({\cal X}(C'))\ar[rrd]^{{\cal Y}_g}
\ar[d]^{{\cal G}(C')({\cal X}(g))}&&
\\{\cal G}(C)({\cal X}(C''))\ar[rr]^{{\cal G}(f)_{{\cal X}(C'')}}
\ar@/_/[rrrr]_{{\cal G}(g\circ f)_{{\cal X}(C'')}}
&&{\cal G}(C')({\cal X}(C''))\ar[rr]^{{\cal G}(g)_{{\cal X}(C'')}}
&&{\cal G}(C'')({\cal X}(C''))
}$$
is commutative because ${\cal G}(f)$ and ${\cal G}(g)$
are natural transformation. The legs of the big triangle compose to
${\cal Y}(g\circ f)$ and the hypothenusis is ${\cal Y}(g)\circ {\cal Y}(f)$.
The equality ${\cal Y}(\id) = \id$ is easy. The contravariant case
is analogous. Q.E.D.
\vskip .08in

Specialize the lemma to the case where $\cal C$ is the (para)cyclic
category. Using the contravariant case of Theorem 2 we get immediately

{\bf Theorem 6.} {\it (assumptions from Theorem 2)
Let $X_\bullet$ be a {\it cyclic} object in ${\cal A}$ with 
boundaries $\partial^i_{X,n}$, degeneracies $\sigma^i_{X,n}$ and
(para)cyclic operators $t_{X,n}$. Let natural transformations $t_n$ from
~(\ref{eq:tcubedeqtcubed}) be denoted by $t^G_n$, with components
$(t^G_n)_M =:t^G_{n,M}$ at $M$ in ${\cal A}$. Then the formulas
\begin{equation}\label{eq:relativecase}\begin{array}{ll}
Y_n : = G^{n+1} X_n, &
\partial_{n,Y}^i 
:= G^i(\epsilon_{G^{n-i}X_n})\circ G^{n+1}(\partial^i_{X,n}),\\
\sigma_{n,Y}^i := G^i(\delta_{G^{n-i}X_n})\circ G^{n+1}(\sigma^i_{X,n}),&
(t_Y)_n := G^{n+1}(t_{X,n}) \circ t^G_{n,X_n},
\end{array}\end{equation}
define a (para)cyclic object $Y_\bullet$ in ${\cal A}$. 
}

The simplicial part of the theorem is 
previously known (cf. e.g.~\cite{Barr:composite}). 
 
{\bf Remark.} 1. Every object $M \in {\cal A}$, gives rise to
a constant simplicial object $X_\bullet$, where $X_n = M$ for
all $n$ and all faces and degeneracies are identities. 
In that case, the assertion of Theorem 2 for $X_\bullet$ is simply the 
main part of Theorem 1. The other parts of the Theorem 1 generalize as well: 
functoriality, augmented version.

2. Monadic version of our results
(giving cosimplicial objects)
is obvious by dualization. We expect that
dihedral etc. analogues of our analysis are possible.

3. For ${\cal A}$ abelian, all flavors of the cyclic homology
associated to the trivial symmetry $\tau^G$ on ${\bf G}$ should not 
carry ``really cyclic information''.
E.g. how do they compare to the 'underlying' cobar homology ?

{\bf Acknowledgements.} 
I thank {\sc M.~Jibladze} for useful conversations and {\sc P.~Bressler} 
for certain motivation (so far still unfullfilled).

{\bf Versions.} We'll address nontrivial examples
and extensions in a later version or a sequel 
to this preprint. This preprint is posted at 
an unusually early stage to facilitate the communication 
with a number of colleagues who expressed 
their interest in the very main construction of this article. 


{\footnotesize
\begin{thebibliography}{99}
\bibitem{Semtriples}
{\sc H.~Appelgate, M.~Barr, J.~Beck et al.},
{\em Seminar on triples and categorical homology theory},
ETH 1966/67, ed. B.~Eckmann, LNM 80, Springer 1969.

\bibitem{Barr:composite} 
{\sc M.~Barr}, {\em Composite triples and derived functors}, 
in~\cite{Semtriples}, pp.~336--356.

\bibitem{Beck:distr} 
{\sc J.~Beck}, {\em Distributive laws}, 
in~\cite{Semtriples}, pp.~119--140.

\bibitem{Bres:cyclic}
{\sc P.~Bressler}, {\em Levels and characters}, a lecture at NOG III,
Mittag-Leffler May 19, 2004; and personal communication on cyclic objects, 
groupoids, gerbes and alia.


\bibitem{BrzNich}
{\sc T.~Brzezi\'{n}ski, F.~F. Nichita}, {\em Yang-Baxter systems
and entwining structures}, {\tt math.QA/0311171}.

\bibitem{bunge:symtopos}
{\sc M.~Bunge, A.~Carboni}, {\em The symmetric topos},
Journal of Pure and Applied Algebra 105 (1995), pp.~233--249.

\bibitem{Cort:DoldKan}
{\sc J.~L. Castiglioni, G.~Corti\~{n}as},
{\em Cosimplicial versus DG-rings: a version of the Dold-Kan correspondence},
J. Pure Appl. Algebra  191  (2004),  no. 1-2, pp.~119--142. 

\bibitem{Connes:book}
{\sc A.~Connes}, {\em Noncommutative geometry}, 
Acad. Press, New York 1994.

\bibitem{DwyerKan:dup}
{\sc W.~G. Dwyer, D.~M. Kan}, {\em Normalizing the cyclic modules of Connes},
Comment. Math. Helv. 60, n.1 (1985), pp.~582--600.

\bibitem{cyclicoper}
{\sc E. Getzler, M. M. Kapranov}, 
{\em Cyclic operads and cyclic homology},
in 'Geometry, topology and physics for Raoul Bott', 
ed. by S-T. Yau, Int. Press 1994.

\bibitem{grandis:symmsimp}
{\sc M.~Grandis}, 
{\em Finite sets and symmetric simplicial sets},
Theory of Applications of Categories (electronic), 
Vol. 8, No. 8, pp.~244--253.

\bibitem{pohu:higherstring}
{\sc Po Hu}, {\em Higher string topology on general spaces},
{\tt math.AT/0401081.}

\bibitem{Loday:cyclicbk}
{\sc J-L.~Loday}, {\em Cyclic homology}, 
Grundl.M.W.~301, Springer 1992, 1998.

\bibitem{LunSko:Hopf} 
{\sc V.~Lunts, Z.~\v{S}koda}, {\em
Hopf modules, flat descent and $Ext$-groups
}, in preparation.

\bibitem{MacLane}
{\sc S.~Mac Lane}, {\em Categories for the working mathematician},
GTM 5, Springer 1971.

\bibitem{Menichi:BV}
{\sc L.~Menichi}, {\it Batalin-Vilkovisky algebras 
and cyclic cohomology of Hopf algebras}, {\tt math.QA/0311276}

\bibitem{Menini:MSRI}
{\sc C.~Menini}, {\em Connections, symmetry operators and 
descent data for triples}, talk at MSRI Hopf Algebras Workshop, 
Oct. 25--28, 1999; video and slides at {\tt www.msri.org}

\bibitem{Nuss}
{\sc P.~Nuss}, {\em Noncommutative descent and non-abelian cohomology},
$K$-Theory  12  (1997),  no. 1, 23--74. 

\bibitem{Skoda:gHGqbun}
{\sc Z.~\v{S}koda}, {\em Globalizing Hopf-Galois extensions} in preparation;
{\em Quantum bundles using coactions and localization}, in preparation.

\bibitem{Skoda:qbun} 
{\sc Z.~\v{S}koda}, {\em Distributive laws for actions of monoidal categories},
{\tt math.CT/0406310}.

\bibitem{Street:monads}
{\sc R.~Street}, {\em The formal theory of monads},
J. Pure Appl. Algebra 2 (1972), pp.~149--168;
\& part II (with {\sc S.~Lack})
J. Pure Appl. Algebra 175 (2002), no. 1-3, pp.~243--265.

\bibitem{Weibel}
{\sc C.~Weibel}, {\em Homological algebra}, 
Cambridge Studies in Adv. Math. 38, Cambridge Univ. Press 1994.

\end{thebibliography}
}
\end{document}